\pdfoutput=1

\RequirePackage{fix-cm}

\documentclass[onecolumn,numbers,sort&compress]{svjour3}

\smartqed

\usepackage{graphicx}
\usepackage{amssymb,amsmath}
\usepackage{txfonts}
\usepackage{natbib}
\usepackage{hyperref}

\hypersetup{pdftitle = The title of my PDF, pdfauthor = My name, pdfsubject= The subject, pdfkeywords = keyword1 keyword2 keyword3}
\hypersetup{colorlinks = true, linkcolor = blue, anchorcolor = red, citecolor = blue, filecolor = red, pagecolor = red, urlcolor = blue}

\journalname{International Journal of Applied and Computational Mathematics}

\begin{document}

\title{Comparing the geometry of the basins of attraction, the speed and the efficiency of several numerical methods}

\author{Euaggelos E. Zotos \and Md Sanam Suraj \and Amit Mittal \and Rajiv Aggarwal}

\institute{Euaggelos E. Zotos: \at
             Department of Physics, School of Science, \\
             Aristotle University of Thessaloniki, \\
             GR-541 24, Thessaloniki, Greece \\
             \email{\url{evzotos@physics.auth.gr}}
         \and
           Md Sanam Suraj: \at
              Department of Mathematics, Sri Aurobindo College, \\
              University of Delhi, Delhi, India \\
%              \email{\url{mdsanamsuraj@gmail.com}}
         \and
           Amit Mittal: \at
              Department of Mathematics, ARSD College, \\
              University of Delhi, New Delhi, India \\
%              \email{\url{to.amitmittal@gmail.com}}
         \and
            Rajiv Aggarwal: \at
              Department of Mathematics, Sri Aurobindo College, \\
              University of Delhi, Delhi, India \\
%              \email{\url{rajiv_agg1973@yahoo.com}}
}

\date{Received: - / Accepted: - / Published online: -}

\titlerunning{Comparing the basins of attraction of several numerical methods}

\authorrunning{Zotos et al.}

\maketitle

\begin{abstract}

We use simple equations in order to compare the basins of attraction on the complex plane, corresponding to a large collection of numerical methods, of several order. Two cases are considered, regarding the total number of the roots, which act as numerical attractors. For both cases we use the iterative schemes for performing a thorough and systematic classification of the nodes on the complex plane. The distributions of the required iterations as well as the probability and their correlations with the corresponding basins of convergence are also discussed. Our numerical calculations suggest that most of the iterative schemes provide relatively similar convergence structures on the complex plane. In addition, several aspects of the numerical methods are compared in an attempt to obtain general conclusions regarding their speed and efficiency. Moreover, we try to determine how the complexity of the each case influences the main characteristics of the numerical methods.


\end{abstract}

\section{Introduction}
\label{intro}

One of the most crucial and fascinating problems in numerical analysis and applied mathematics is solving an nonlinear equation. A large number of research papers is available for computing and solving one nonlinear equation or even a system of nonlinear equations, using numerical methods (see e.g., \cite{CNK14,CN150,CN17a,CN17b,CN16,CN18,NSC12}). Numerical solutions of systems of nonlinear equations are also applicable in chaotic dynamics (see e.g., \cite{HMB16,MHBA15}). There are two basic techniques for finding the roots of the equation, that is the direct method and the iterative method. Moreover, the iterative method converges to a root of the particular function, under certain initial conditions. However, the rate of convergence highly depends on the particular numerical technique used. The iterative technique which exhibits very rapid convergence, requires stronger conditions to guarantee convergence. Therefore, one can  also classify and compare the iterative methods by their order of convergence, as well as by their complexity, expressed by the number of function evaluated per step.

There are various ways in which we can compare the iterative techniques, proposed for solving nonlinear equations. The main problem with the iterative methods is the fact that these methods may not converge or may converge to another root of the equation. However, in general terms, if a good starting point is given then this problem will not arise. In addition, the initial point is chosen randomly and this can be improved by choosing a number other than the randomly chosen initial point which provides an incrementally better method of comparison. This motivates us to develop a better and effective methodology for the comparison of algorithms for solving nonlinear equations. The methods are assorted on the basis of their order, informational efficiency as well as efficiency index. We also consider another measures by discussing the basins of attraction of the method and its dependence on the order.

The basins of attraction, associated with the roots of an equation, provide qualitative information of the intrinsic properties of the equation. It was \cite{S01} who introduced the idea of the basins of attraction of root finding methods. Furthermore, a comparison between the classical Newton's method and the third order methods described in \cite{H64} and \cite{P80}, was provided. The aim of the present paper is to try to give some answers to the following questions: (i) how can we compare the basins of attraction of different algorithms with the same order of convergence? (ii) how can we compare the basins of attraction of different algorithms with different order of convergence? (iii) Are there any relationships connecting the order of the numerical methods with their corresponding speed and efficiency?

Recently, \cite{SNC11} discussed ten numerical methods of various orders, in which they included optimal methods of order two, four, eight and sixteen, where the method of order $2^n$ is optimal in the sense that it requires $n+1$ function and derivative evaluations per cycle. They also discussed the basins of attraction and shown that the basins corresponding to the Kung-Traub method and all methods based on it have more chaotic and more fractal basin boundaries than those of King's method and the methods based on it. \cite{CN15a} studied the basins of attraction corresponding to the Euler-Cauchy's method, using seven different polynomials having multiple roots, with fixed multiplicity. In the same vein, \cite{CN15b} compared the basins of attraction of Kanwar-Bhatia-Kansal family in order to determine the best numerical method of fourth order.

In our study, we consider sixteen numerical methods of different order (varying between 2 and 16), where ten of them were taken from \cite{SNC11}. We use these iterative formulae for revealing the basins of attraction on the complex plane, corresponding to the roots of simple polynomial equations. Moreover we are going to use several quantitative information for comparing the speed of convergence as well as the efficiency of these numerical methods. On this basis, we may argue that the presented numerical outcomes are novel and this is exactly the contribution of our work.

The paper is structured as follows: The main numerical results, regarding the basins of attraction of the several numerical methods, are presented in Section \ref{num}. In Section \ref{base} we use the basin entropy for comparing the degree of the fractality and the chaoticity of the fractal basin boundaries. Our article ends with Section \ref{conc}, where the most important conclusions of our numerical analysis are emphasized.

\section{Numerical methods}
\label{num}

In the literature one can find a plethora of methods for numerically solving an equation with only one variable parameter. In this study we will consider and therefore compare sixteen methods, whose order of convergence is varying from 2 to 16. The methods under consideration are the following:
\begin{enumerate}
  \item The Newton-Raphson's optimal method of second order \cite{CdB73}.
  \item The Halley's method of third order \cite{H64}.
  \item The Chebyshev's method of third order \cite{T64}
  \item The super Halley's method of fourth order \cite{GH01}.
  \item The modified super Halley's optimal method of fourth order \cite{CH08}.
  \item The King's method of fourth order \cite{K73}.
  \item The Jarratt's method of fourth order \cite{J66}.
  \item The Kung-Traub's optimal method of fourth order \cite{KT74}.
  \item The Maheshwari's optimal method of fourth order \cite{M09}.
  \item The Murakami's method of fifth order \cite{M78}.
  \item The Neta's method of sixth order \cite{N79}.
  \item The Chun-Neta's method of sixth order \cite{CN12}.
  \item The Neta-Johnson's method of eighth order \cite{NJ08}.
  \item The Neta-Petkovic's optimal method of eighth order \cite{NP10}.
  \item The Neta's method of fourteenth order \cite{N81}.
  \item The Neta's method of sixteenth order \cite{N81}.
\end{enumerate}
The analytical expressions of all the above-mentioned iterative schemes can be found in the Appendix of \cite{Z17}.

\begin{figure*}
\centering
\resizebox{\hsize}{!}{\includegraphics{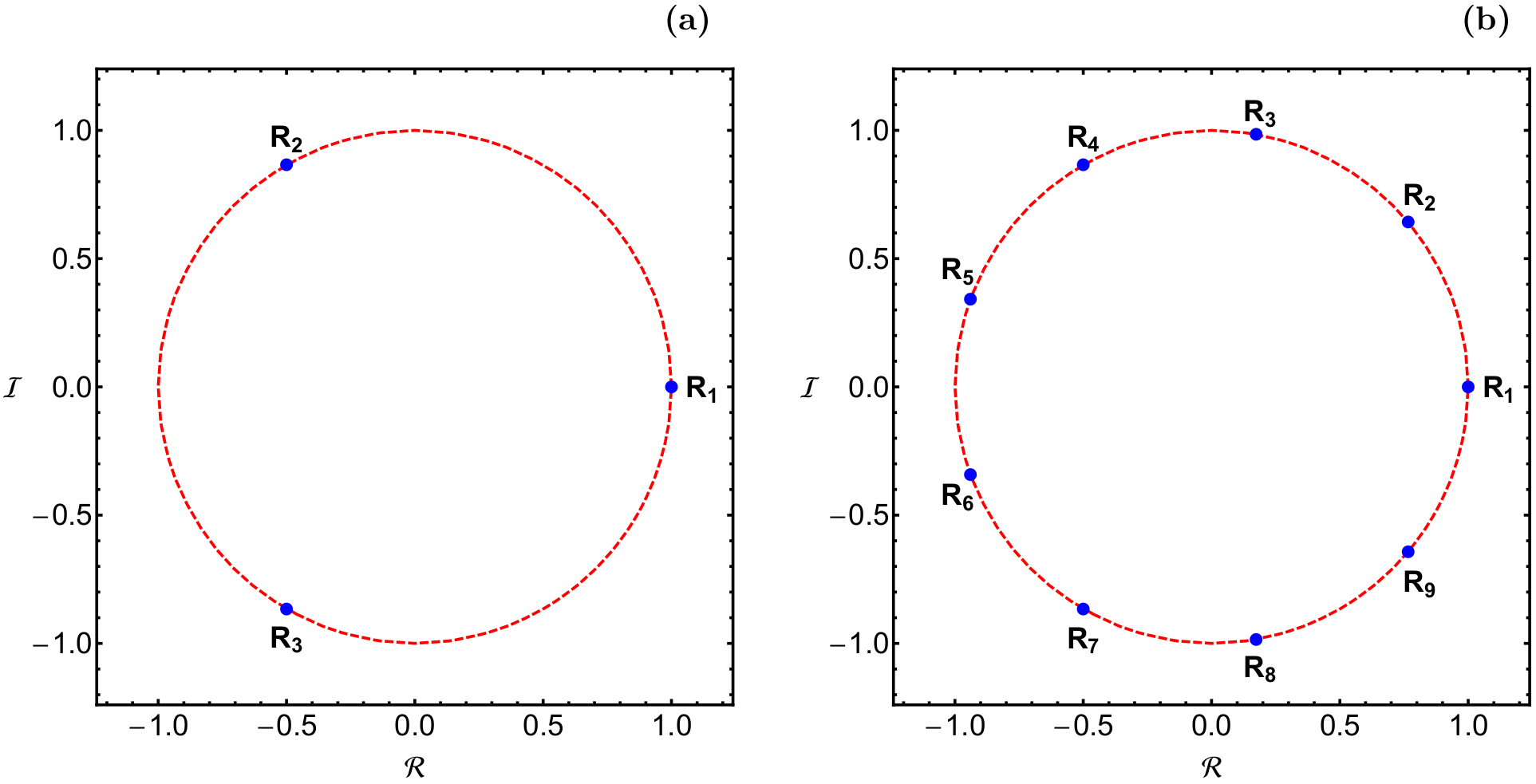}}
\caption{The position of the roots on the complex plane for the equation (a-left): $f_1(z) = z^3 - 1 = 0$ and (b-right): $f_2(z) = z^9 - 1 = 0$. In both cases, all the roots lie on an circle with radius $R = 1$, which is indicated by red dashed line. (Color figure online).}
\label{roots}
\end{figure*}

For comparing the convergence properties as well as the efficiency of the numerical methods we will use them in two different cases, regarding the total number of the existing roots (numerical attractors). In particular, we will numerically solve the simple polynomial functions $f_1(z) = z^3 - 1 = 0$ and $f_2(z) = z^9 -  1 = 0$. The first function $f_1(z)$ has three roots (one real and two complex conjugate), while the second function $f_2(z)$ has nine roots (one real and eight complex conjugate). In both cases, all the roots lie on an circle with radius $R = 1$ (see Fig. \ref{roots}). Furthermore, it is interesting to note that both complex equations have three roots in common, that is $z = 1$ and $z = - 1/2 \pm i \sqrt{3}/2$.

The philosophy behind all the numerical method is the following: An initial complex number $z = \mathcal{R} + i \mathcal{I}$ on the complex plane activates the code, while the iterative procedure continues until a root (attractor) is reached, with the desired predefined accuracy. If the particular initial condition leads to one of the roots it means that the numerical method converges for that particular initial condition $(\mathcal{R}, \mathcal{I})$. At this point, it should be emphasized that in general terms a method does not converge equally well for all the available initial conditions. The sets of the initial conditions which lead to the same final state (root) compose the so-called basins of convergence/attraction or convergence/attracting domains/regions. Nevertheless, it should be clarified that the basins of attraction, associated with the roots of the polynomials, should not be mistaken, by no means, with the basins of attraction which are present in dissipative systems.

A double scan of the complex plane is performed for revealing the structures of the basins of convergence. In particular, a dense uniform grid of $1024 \times 1024$ $(\mathcal{R},\mathcal{I})$ nodes is defined which shall be used as initial conditions of the iterative schemes. The number $N$ of the iterations, required for obtaining the desired accuracy, is also monitored during the classification of the nodes. For our computations, the maximum allowed number of iterations is $N_{\rm max} = 500$. Moreover the iterations stop only when a root is reached, with accuracy of $10^{-15}$ for both real and imaginary parts.

In the following subsections we will compare the basins of attraction of the several numerical methods, by considering two cases regarding the total nature of the roots. For the classification of the nodes on the complex plane we will use color-coded diagrams (CCDs), in which each pixel is assigned a different color, according to the final state (root) of the corresponding initial condition. These color-coded diagrams are also known in the bibliography as ``polynomiographs". It should be noted that all the computational tools as well as the numerical approaches are as in \cite{Z17,Z18}.

\begin{figure*}[!t]
\centering
\resizebox{\hsize}{!}{\includegraphics{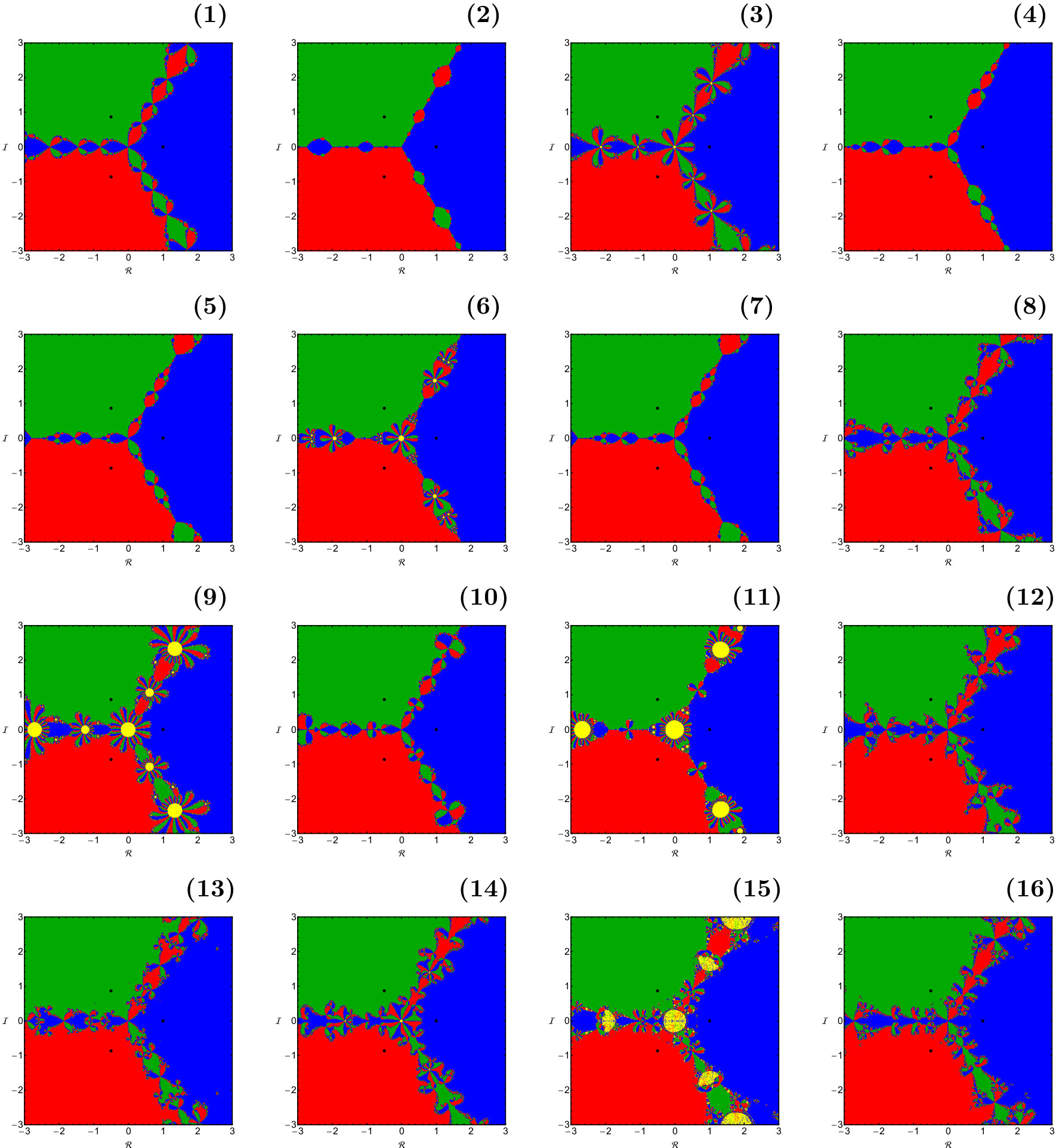}}
\caption{The basins of attraction on the complex plane for the first case, where three roots exist. The positions of the three roots are indicated by black dots. The color code is as follows: $R_1$ root (blue); $R_2$ root (green); $R_3$ root (red); initial conditions for which the iterative schemes lead to extremely large numbers (yellow); initial conditions for which the iterative schemes immediately abort (orange). The numbers of the panels correspond to the numerical formulae, as they have been listed at the beginning of Section \ref{num}. (Color figure online).}
\label{c1}
\end{figure*}

\begin{figure*}[!t]
\centering
\resizebox{\hsize}{!}{\includegraphics{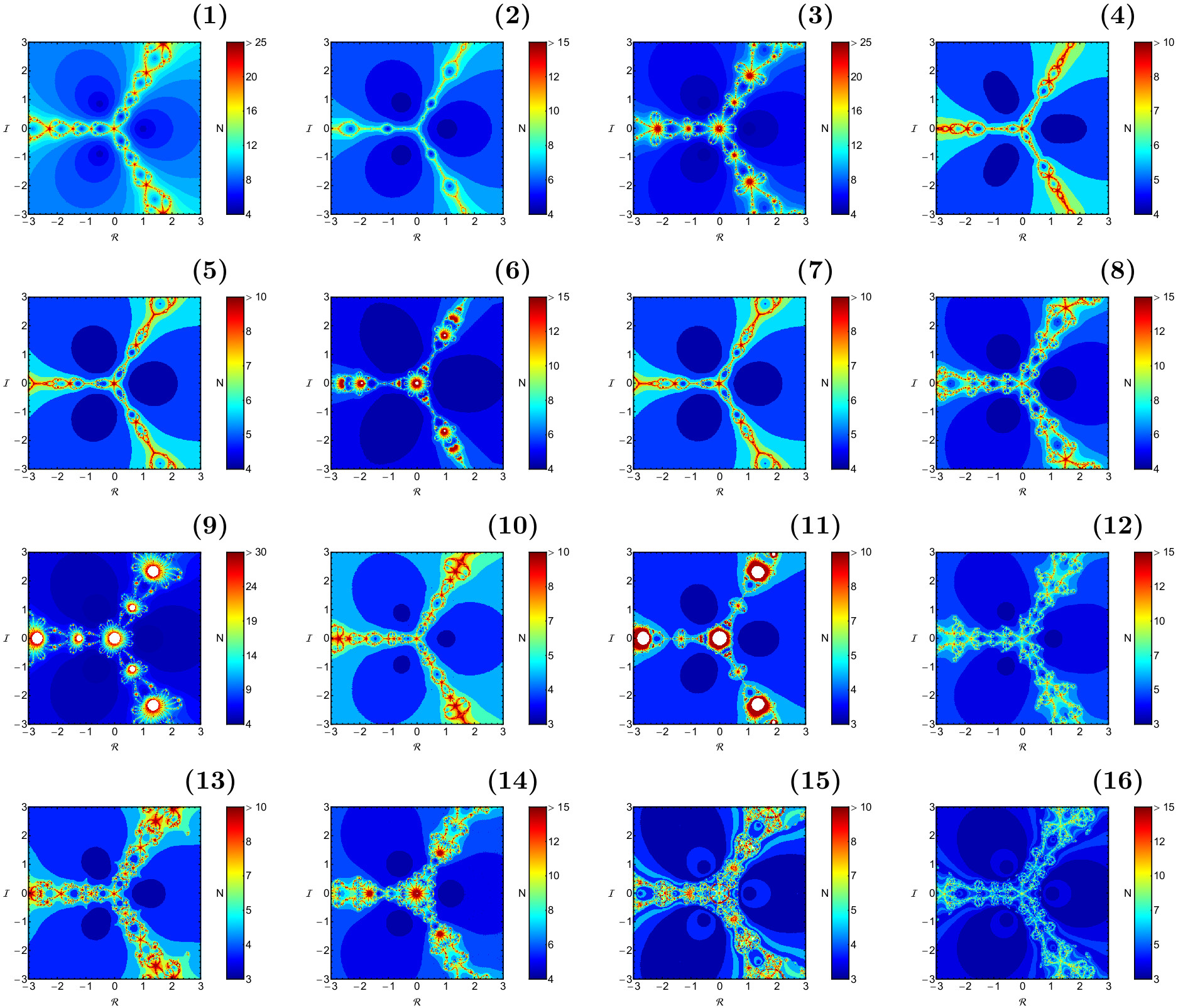}}
\caption{The corresponding distributions of the number $N$ of the required iterations for obtaining the basins of attraction shown in Fig. \ref{c1}. All ill-behaved initial conditions are shown in white. (Color figure online).}
\label{n1}
\end{figure*}

\begin{figure*}[!t]
\centering
\resizebox{\hsize}{!}{\includegraphics{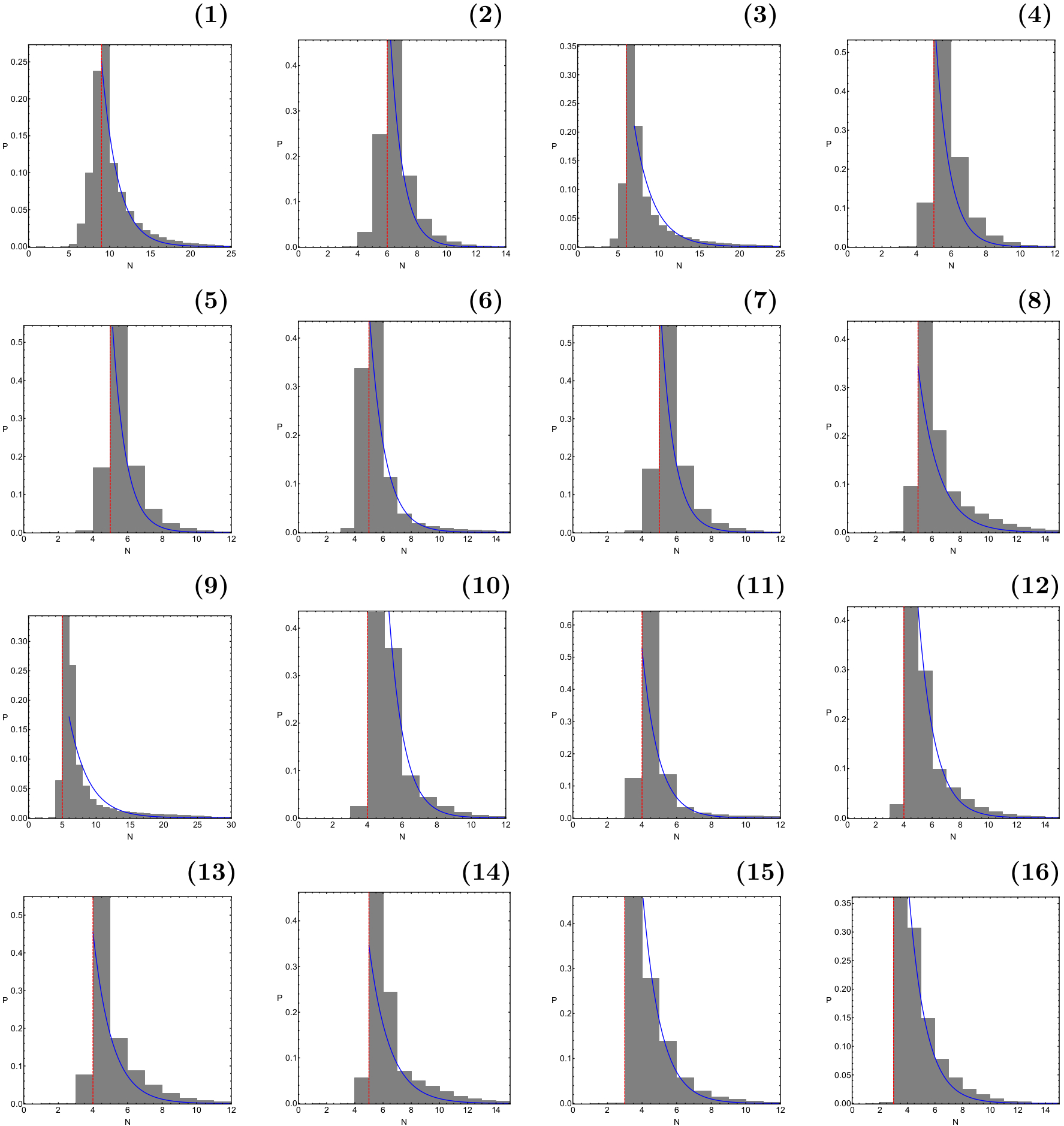}}
\caption{The corresponding probability distributions of the required iterations for obtaining the basins of attraction shown in Fig. \ref{c1}. The vertical, dashed, red lines indicate, in each case, the most probable number $N^{*}$ of iterations. (Color figure online).}
\label{p1}
\end{figure*}

\begin{figure*}[!t]
\centering
\resizebox{\hsize}{!}{\includegraphics{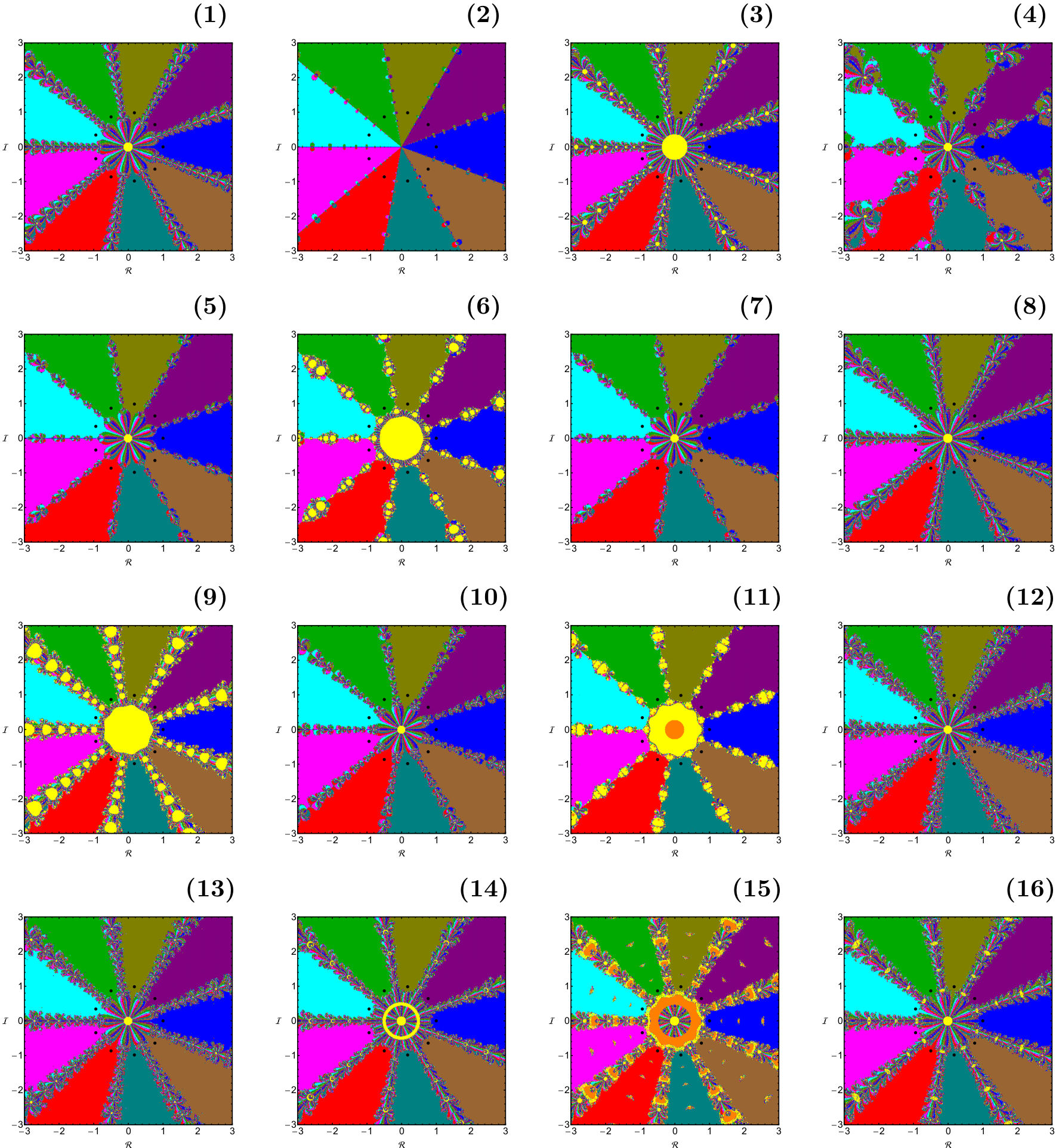}}
\caption{The basins of attraction on the complex plane for the first case, where nine roots exist. The positions of the three roots are indicated by black dots. The color code is as follows: $R_1$ root (blue); $R_2$ root (purple); $R_3$ root (olive); $R_4$ root (green); $R_5$ root (cyan); $R_6$ root (magenta); $R_7$ root (red); $R_8$ root (teal); $R_9$ root (brown); initial conditions for which the iterative schemes lead to extremely large numbers (yellow); initial conditions for which the iterative schemes immediately abort (orange). The numbers of the panels correspond to the numerical formulae, as they have been listed at the beginning of Section \ref{num}. (Color figure online).}
\label{c2}
\end{figure*}

\begin{figure*}[!t]
\centering
\resizebox{\hsize}{!}{\includegraphics{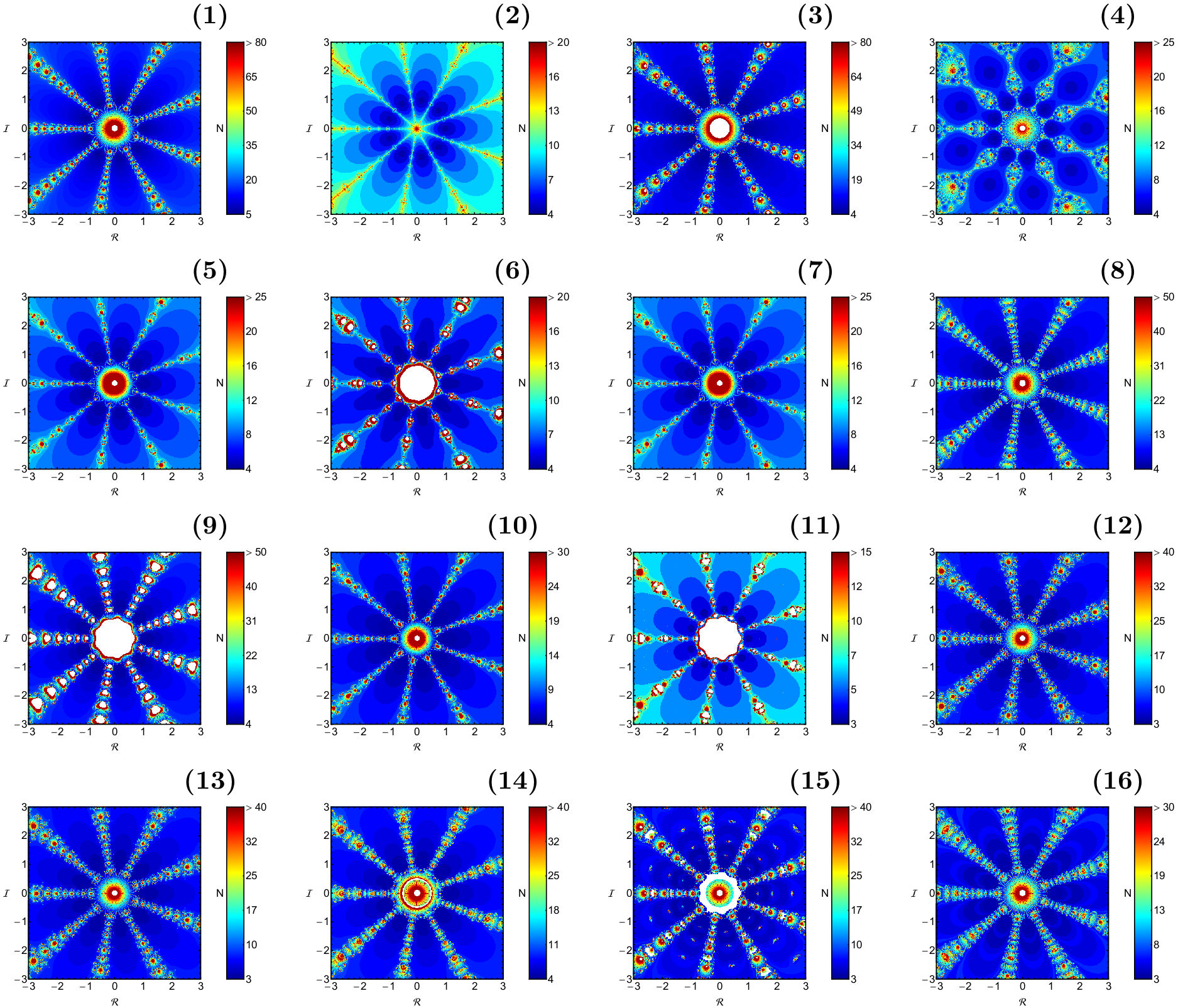}}
\caption{The corresponding distributions of the number $N$ of the required iterations for obtaining the basins of attraction shown in Fig. \ref{c2}. All ill-behaved initial conditions are shown in white. (Color figure online).}
\label{n2}
\end{figure*}

\begin{figure*}[!t]
\centering
\resizebox{\hsize}{!}{\includegraphics{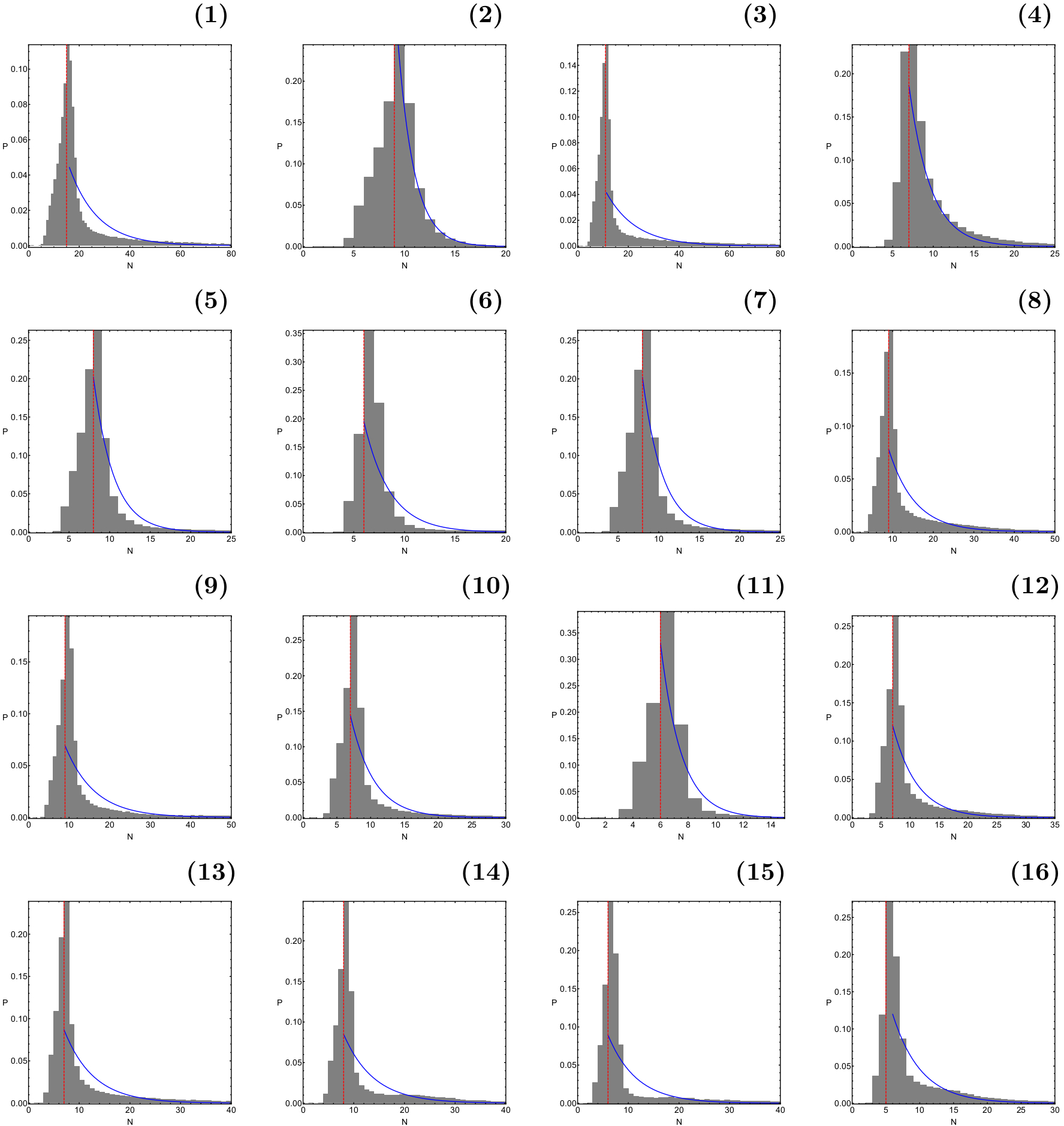}}
\caption{The corresponding probability distributions of the required iterations for obtaining the basins of attraction shown in Fig. \ref{c2}. The vertical, dashed, red lines indicate, in each case, the most probable number $N^{*}$ of iterations. (Color figure online).}
\label{p2}
\end{figure*}

\subsection{Case I: Three roots (numerical attractors)}
\label{ss1}

We begin with the first case, where the complex equation $f_1(z)$ has three roots. The basins of attraction on the complex plane, corresponding to the sixteen iterative schemes are presented in Fig. \ref{c1}. It is seen that, in general terms, the geometry of the convergence regions around the roots is almost identical to all the numerical methods. On the other hand, the geometry of the regions, mainly in the vicinity of the basin boundaries, is completely different and highly fractal\footnote{It should be noted that when we state that an area is fractal we simply mean that it displays a fractal-like geometry, without using (at least for now) any quantitative arguments, such as the fractal dimension, as in \cite{AVS01,AVS09}.} This directly implies that for the initial conditions inside these chaotic domains it is extremely difficult (or even impossible) to know beforehand their final state (root). In all cases the basins of attraction, associated with the three roots, extend to infinity.

For several numerical methods (i.e., Chebyshev, King, Maheshwari, Neta, Neta-Petkovic, Neta method of fourteenth order) we detected a non-zero amount of initial conditions for which the corresponding iterative schemes lead quickly to extremely large numbers. Additional computations revealed that for these initial conditions in some cases the iterative formulae, after a considerable amount of iterations, do eventually converge to one of the roots. However, for the majority of these initial conditions the iterative formulae do not converge to any of the roots. It is seen in Fig. \ref{c1} that these rebel initial conditions are situated in the areas of the complex plane with the most complicated geometry. In particular, they form additional spherical basins (yellow regions), which are mainly located at the centers of the flower-shaped basins of attraction.

It should be noted that in the case of the Neta method of fourteenth order (see panel (15) of Fig. \ref{c1}) we encountered a small portion of initial conditions for which the iterative method aborts at the first iteration.

In Fig. \ref{n1} we provide the corresponding distributions of the number $N$ of the required iterations for obtaining the basins of attraction shown in Fig. \ref{c1}. In almost all cases we observe the expected behavior according to which the fastest converging nodes are those with initial conditions inside the basins of attraction, while the slowest initial conditions are those located in the vicinity of the basin boundaries.

The corresponding probability distributions of iterations are given in Fig. \ref{p1}. The definition of the probability $P$ is the following: assume that after $N$ iterations $N_0$ initial conditions on the complex plane converge to one of the roots of the system. Then $P = N_0/N_t$, where $N_t$ is the total number of initial conditions in every CCD. Our results suggest that for almost all the numerical methods more than 98\% of the initial conditions converge within the first 20 iterations, while only for the Newton-Raphson and the Maheshwari methods the tails of the corresponding histograms extends to more than 20 iterations. The most probable number of iterations $N^{*}$ (see the vertical, dashed, red lines in the histograms) seems, in general terms, to decrease as we proceed to numerical methods of higher order.

\subsection{Case II: Nine roots (numerical attractors)}
\label{ss2}

In this case, the complex equation $f_2(z)$ has nine roots and therefore the convergence structure of the complex plane should be much more complicated. In Fig. \ref{c2} we provide the basin diagrams for the sixteen numerical methods. As expected the geometry of basins of attraction is more complicated with respect to what observed in Fig. \ref{c1}, regarding the case with only three roots. In particular, in most of the cases the central region of the CCDs is occupied by a highly fractal mixture of initial conditions with a shape that resembles a flower with many petals.
Looking at panels of Fig. \ref{c2} it becomes evident that there are two main differences, regarding the ill-behaved initial conditions. More specifically:
\begin{itemize}
  \item Initial conditions for which the iterative schemes lead, at least temporarily, at extremely large numbers are present in almost all cases. Indeed such initial conditions were identified for numerical methods apart from the Halley method. Once more, there initial conditions form circular basins, mainly at the intersections of the usual basins of attraction. The most prominent basins of this type of initial conditions are observed for the following types of numerical methods: Chebyshev, King, Maheshwari, Neta, Neta-Petkovic and Neta method of fourteenth order. In all other cases, these initial conditions are mainly located near the origin.
  \item For the numerical methods King, Maheshwari, Neta, Neta-Petkovic and Neta of fourteenth order we identified initial conditions for which the corresponding numerical formulae abort at the first step of the iterative procedure. In fact for the Neta and Neta of fourteenth order methods it is seen in panels (11) and (15) of Fig. \ref{c2}, respectively that the corresponding initial conditions form circular basins near the origin.
\end{itemize}

The corresponding distributions of the number of iterations and the probability are given in Figs. \ref{n2} and \ref{p2}, respectively. In Fig. \ref{p2} one can see that in this case (where nine roots exist) the tails of almost all histograms extend to higher values than those of the previous case (see Fig. \ref{p1}). Only the histograms corresponding to Halley and the super Halley methods display a smooth distribution, without the presence of long extended tails.

\begin{figure*}[!t]
\centering
\resizebox{\hsize}{!}{\includegraphics{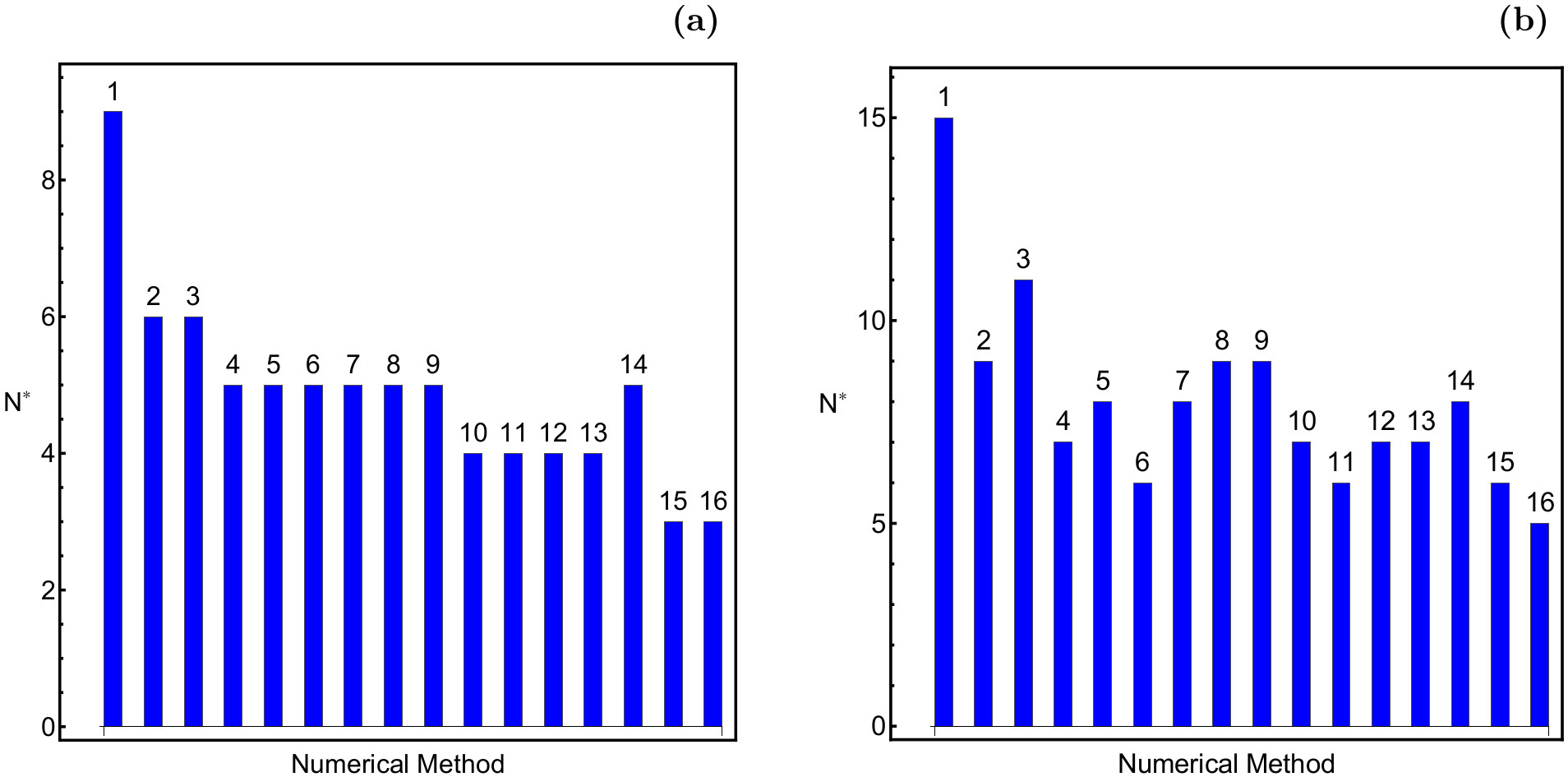}}
\caption{The distribution of the most probable number of iterations $N^{*}$ for all the numerical methods for the case where (a-left): three roots exist and (b-right): nine roots exist. The numbers of each bar correspond to the numerical formulae, as they have been listed at the beginning of Section \ref{num}. (Color figure online).}
\label{iter}
\end{figure*}

\subsection{Comparison of the numerical methods}
\label{comp}

In this subsection we will compare several aspects of the different numerical methods, derived from the two cases presented earlier in subsections \ref{ss1} and \ref{ss2}.

One of the most important quantities is the most probable number of iterations $N^{*}$. In Fig. \ref{iter}(a-b) we present the corresponding values of $N^{*}$ for all the numerical methods. In panel (a), which corresponds to the case of three roots, we see that $N^{*}$ displays a clear reduction as we proceed to numerical methods of higher order. In other words, the speed of the convergence increases, since less iterations are required for obtaining the same accuracy. Moreover, it should be noted that for most of the numerical methods of the same order the value of $N^{*}$ is the same. On the other hand, in the case where nine roots exist (see panel (b) of Fig. \ref{iter}) things are not so clear, regarding the speed of the convergence. However, we may argue again that, in general terms, the value of $N^{*}$ decreases with increasing order. For example, for the Newton-Raphson method of second order we have that $N^{*} = 15$, while for the Neta method of sixteenth order $N^{*} = 5$.

The probability distributions, shown in Figs. \ref{p1} and \ref{p2}, can provide additional useful information. For instance, obtaining the best fit of the tails of the histograms would give us answers, regarding the convergence of the numerical methods. The Laplace distribution is used for fitting the tails of the histograms. The choice of the particular type of distribution is the most natural choice, if we take into consideration that the Laplace distribution is very common in system where transient chaos occurs (e.g., \cite{ML01,SASL06,SS08}).

The function of the probability density (PDF) of the Laplace distribution can be computed as follows
\begin{equation}
P(N | a,b) = \frac{1}{2b}
 \begin{cases}
      \exp\left(- \frac{a - N}{b} \right), & \text{if } N < a \\
      \exp\left(- \frac{N - a}{b} \right), & \text{if } N \geq a
 \end{cases},
\label{pdf}
\end{equation}
where $a$ is the so-called location parameter, while $b > 0$, is the usually known as the diversity. Since the tails of the histograms refer to $N > N^{*}$ in our case we are only interested for the $x \geq a$ part of the PDF. The values of the location parameter as well as the diversity for both cases, regarding the number of roots, are given in Table \ref{tab1}.

\begin{table}[!ht]
\begin{center}
   \caption{The most probable number of iterations $N^{*}$, the location parameter $a$ and the diversity $b$, for all the cases presented in Figs. \ref{p1} and \ref{p2}.}
   \label{tab1}
   \setlength{\tabcolsep}{10pt}
   \begin{tabular}{@{}lccc}
      \hline
      Figure & $N^{*}$ & $a$ & $b$ \\
      \hline
      \ref{p1} Panel  (1) & 9 &     $N^{*}$ & 1.99 \\
      \ref{p1} Panel  (2) & 6 &     $N^{*}$ & 0.83 \\
      \ref{p1} Panel  (3) & 6 & $N^{*} + 1$ & 2.37 \\
      \ref{p1} Panel  (4) & 5 &     $N^{*}$ & 0.81 \\
      \ref{p1} Panel  (5) & 5 &     $N^{*}$ & 0.77 \\
      \ref{p1} Panel  (6) & 5 &     $N^{*}$ & 1.08 \\
      \ref{p1} Panel  (7) & 5 &     $N^{*}$ & 0.77 \\
      \ref{p1} Panel  (8) & 5 &     $N^{*}$ & 1.45 \\
      \ref{p1} Panel  (9) & 5 & $N^{*} + 1$ & 2.91 \\
      \ref{p1} Panel (10) & 4 & $N^{*} + 1$ & 0.87 \\
      \ref{p1} Panel (11) & 4 &     $N^{*}$ & 0.94 \\
      \ref{p1} Panel (12) & 4 & $N^{*} + 1$ & 1.69 \\
      \ref{p1} Panel (13) & 4 &     $N^{*}$ & 1.09 \\
      \ref{p1} Panel (14) & 5 &     $N^{*}$ & 1.44 \\
      \ref{p1} Panel (15) & 3 & $N^{*} + 1$ & 1.03 \\
      \ref{p1} Panel (16) & 3 & $N^{*} + 1$ & 1.24 \\
      \hline
      \hline
      \ref{p2} Panel  (1) & 15 & $N^{*} + 1$ & 11.14 \\
      \ref{p2} Panel  (2) &  9 &     $N^{*}$ &  1.61 \\
      \ref{p2} Panel  (3) & 11 &     $N^{*}$ & 11.70 \\
      \ref{p2} Panel  (4) &  7 &     $N^{*}$ &  2.68 \\
      \ref{p2} Panel  (5) &  8 &     $N^{*}$ &  2.47 \\
      \ref{p2} Panel  (6) &  6 &     $N^{*}$ &  2.56 \\
      \ref{p2} Panel  (7) &  8 &     $N^{*}$ &  2.47 \\
      \ref{p2} Panel  (8) &  9 &     $N^{*}$ &  6.41 \\
      \ref{p2} Panel  (9) &  9 &     $N^{*}$ &  7.22 \\
      \ref{p2} Panel (10) &  7 &     $N^{*}$ &  3.47 \\
      \ref{p2} Panel (11) &  6 &     $N^{*}$ &  1.51 \\
      \ref{p2} Panel (12) &  7 &     $N^{*}$ &  4.16 \\
      \ref{p2} Panel (13) &  7 &     $N^{*}$ &  5.75 \\
      \ref{p2} Panel (14) &  8 &     $N^{*}$ &  5.88 \\
      \ref{p2} Panel (15) &  6 &     $N^{*}$ &  5.60 \\
      \ref{p2} Panel (16) &  5 & $N^{*} + 1$ &  4.18 \\
      \hline
   \end{tabular}
\end{center}
\end{table}

\begin{figure*}[!t]
\centering
\resizebox{\hsize}{!}{\includegraphics{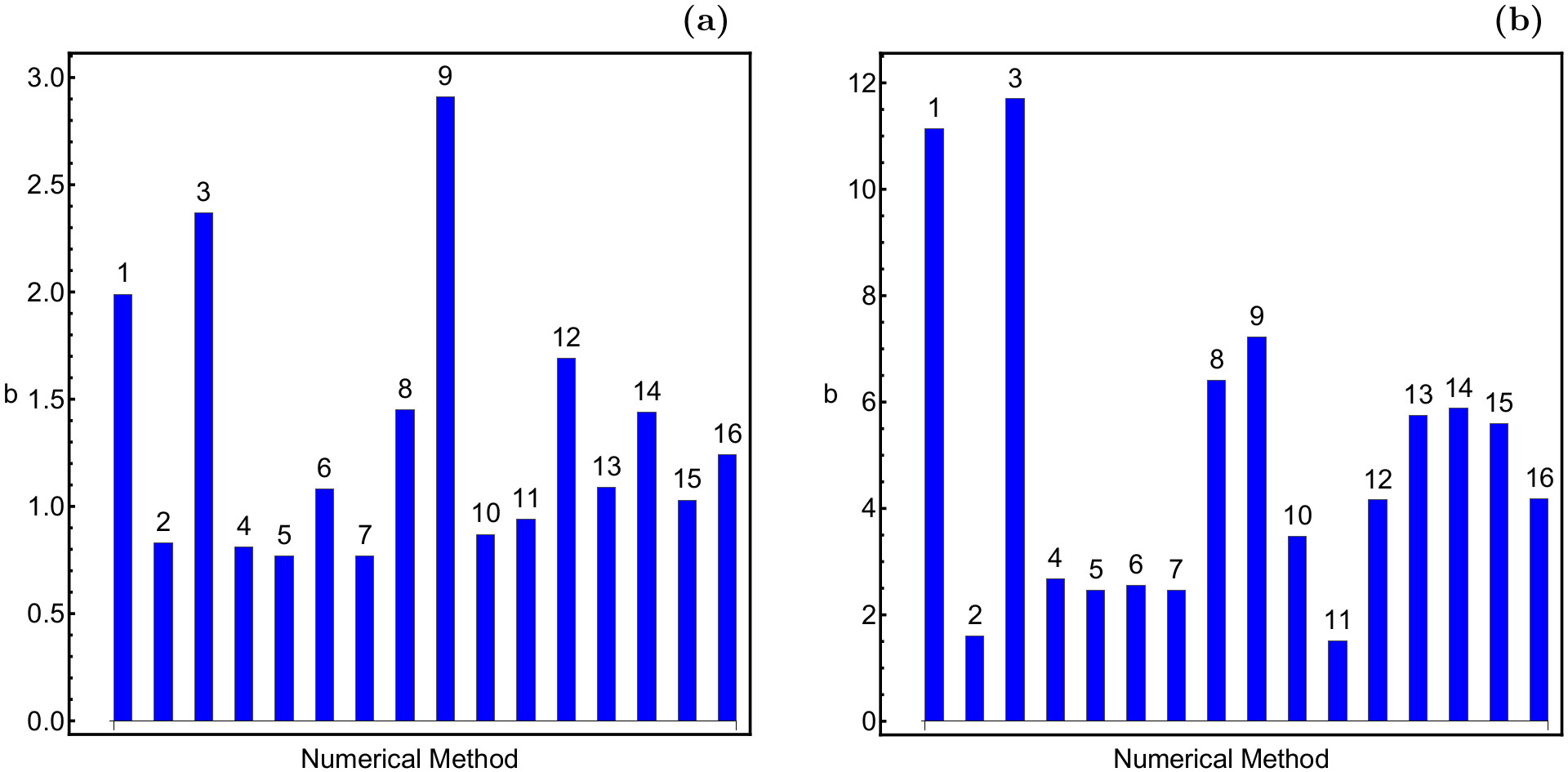}}
\caption{The distribution of the diversity $b$ of the Laplace PDF for all the numerical methods for the case where (a-left): three roots exist and (b-right): nine roots exist. The numbers of each bar correspond to the numerical formulae, as they have been listed at the beginning of Section \ref{num}. (Color figure online).}
\label{div}
\end{figure*}

\begin{figure*}[!t]
\centering
\resizebox{\hsize}{!}{\includegraphics{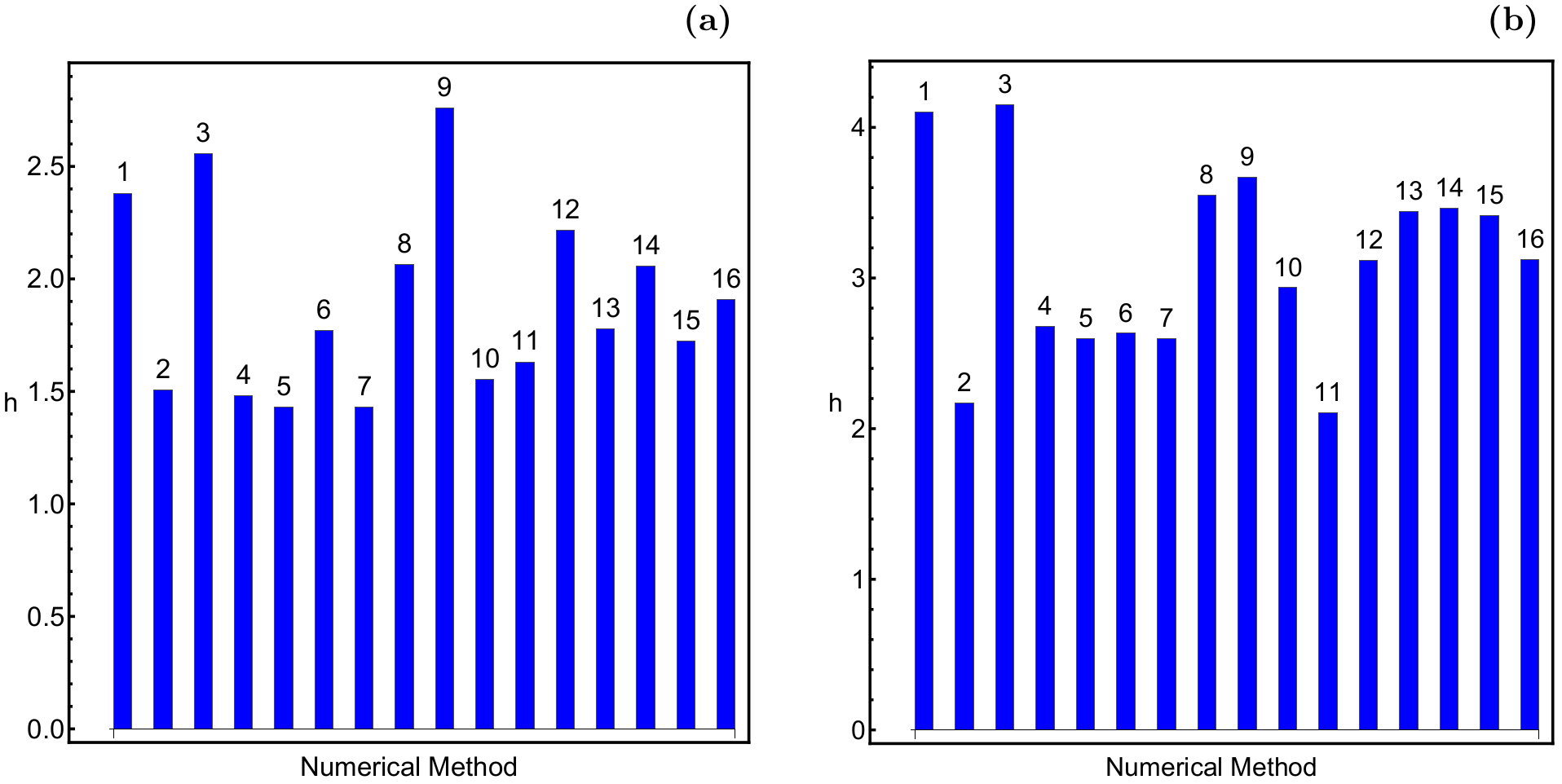}}
\caption{The distribution of the differential entropy $h$ of the Laplace PDF for all the numerical methods for the case where (a-left): three roots exist and (b-right): nine roots exist. The numbers of each bar correspond to the numerical formulae, as they have been listed at the beginning of Section \ref{num}. (Color figure online).}
\label{deh}
\end{figure*}

It is observed, that for most of the cases the number $N^{*}$ is very close to the location parameter $a$. Especially when nine roots exist these two quantities completely coincide for the majority of the numerical methods. In Fig. \ref{div} we present the distribution of the diversity $b$ for both studied cases. We observe in panel (a), which correspond to the first case where three roots exist, that for almost all the numerical methods $b$ is lower than 2, while on the other hand in panel (b), which correspond to the second case where nine roots exist, for almost all the iterative schemes $b > 2$. Here, it should be noted that a low value of the diversity $(b < 2)$ implies that the corresponding distribution of the histogram is compact and well organized around the maximum value $N^{*}$. On the other hand, when $b > 2$ the corresponding distribution is much more extended and therefore occupies a larger range of $N$ values. This is exactly the reason of why the Laplace PDF fails to properly fit almost all cases in Fig. \ref{p2}, while the same type of distribution is the best fit for almost all cases shown in Fig. \ref{p1}.

The diversity is directly related with another important quantity which is called the ``differential entropy". The differential entropy measures, to a continuous probability distribution, the average surprisal of a random variable. In the case of the Laplace probability distribution the differential entropy is defined as $h = 1 + \ln(2b)$. In panels (a-b) of Fig. \ref{deh} we illustrate the values of the differential entropy for all the numerical methods and for both cases. It is seen, that when three roots exist the highest value of the differential entropy is measured for the Maheshwari's optimal method, while for the case of nine roots the highest value of $h$ corresponds to the Chebyshev's method. Obviously, the highest as well as lowest values of $h$ are for those methods where the diversity $b$ displays the highest and lowest values, respectively.

Before closing this section we would like to shed some light to the issue of ill-behaved initial conditions. In the previous two subsections we mentioned the presence of initial conditions for which the numerical code aborts at the very first step of the iterative procedure. In order to demystify the reason of this malfunction we chose several of these ill-behaved initial conditions and we monitored all the computed quantities. It was found that for these initial conditions several terms, entering the iterative formulae, tend either to extremely large or to extremely small numbers. This automatically implies that there are cases in which a division by zero (practically an extremely small number, of the order of $10^{-16}$, which computationally is equal to zero) is encountered and therefore this is the reason of why the numerical code immediately aborts.

\section{The basin entropy}
\label{base}

So far, in the numerical results presented in the previous Section, we used only qualitative arguments for discussing the degree of the fractality of the basins of convergence on the complex plane. However, there is no doubt that quantitative results, regarding the evolution of the fractality, would be very informative. Very recently, a new quantitative tool was introduced, for measuring the degree of the basin fractality \cite{DWGGS16,DWGGS18}. This new dynamical quantity is called ``basin entropy" and it measures the degree of fractality (or unpredictability) of the basins, by examining their topological properties.

The basin entropy works according to the following numerical algorithm. If there are $N(A)$ attractors (roots) in a certain region $R = [-3,3] \times [-3,3]$ on the complex plane, then we subdivide $R$ into a grid of $N$ square boxes, where each cell of the gird may contain between 1 and $N(A)$ attractors. Then the probability that inside the cell $i$ the corresponding attractor is $j$ is denoted by $P_{i,j}$. Taking into account that inside each cell the initial conditions are completely independent, the Gibbs entropy, of every cell $i$ reads
\begin{equation}
S_{i} = \sum_{j=1}^{m_{i}}P_{i,j}\log_{10}\left(\frac{1}{P_{i,j}}\right),
\end{equation}
where $m_{i} \in [1,N_{A}]$ is the total number of the attractors inside the cell $i$.

The total entropy of the entire region $R$, on the complex plane, can easily be calculated by adding the entropies of the $N$ cells of the grid as $S = \sum_{i=1}^{N} S_{i}$. Therefore, the total entropy, corresponding to the total number of cells $N$ is called basin entropy and it is given by
\begin{equation}
S_{b} = \frac{1}{N}\sum_{i=1}^{N}\sum_{j=1}^{m_{i}}P_{i,j}\log_{10}\left(\frac{1}{P_{i,j}}\right).
\end{equation}

\begin{figure*}[!t]
\centering
\resizebox{\hsize}{!}{\includegraphics{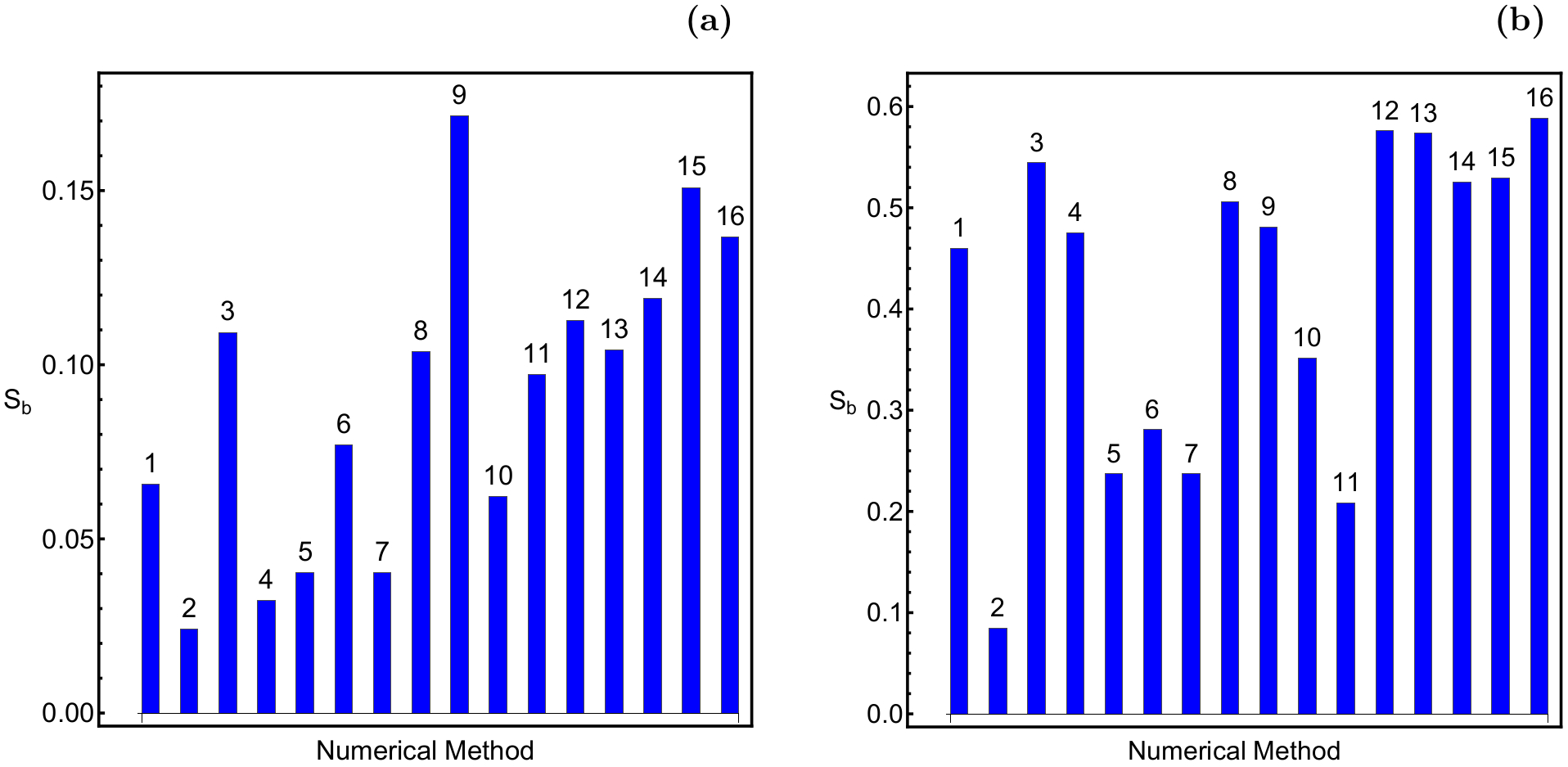}}
\caption{The distribution of the values of the basin entropy $S_b$ of the complex plane, for all the numerical methods for the case where (a-left): three roots exist and (b-right): nine roots exist. The numbers of each bar correspond to the numerical formulae, as they have been listed at the beginning of Section \ref{num}. (Color figure online).}
\label{be}
\end{figure*}

Following the above-mentioned algorithm and also using the value $\varepsilon = 0.005$, suggested in \cite{DWGGS16}, we calculated the numerical value of the basin entropy $S_b$ of the complex plane, for all the numerical methods. At this point, it should be emphasized that the initial conditions, on the complex plane, for which the iterative schemes fail to converge to one of the roots, were counted as an additional type of basin, which coexist along with the regular basins of attraction.
The numerical values of $S_b$, corresponding to the sixteen iterative schemes, are given in the bar chart of Fig. \ref{be}. In panel (a), which corresponds to the case of three roots, one can see that the lowest value of $S_b$ was measured for the Halley method, while the largest value of the basin entropy corresponds to the Maheshwari's method. For the case where nine roots are present (see panel (b) of Fig. \ref{be}) the lowest value of $S_b$ corresponds again to the Halley method, while the convergence planes with the highest degree of fractality are those produced using the methods: Chun-Neta, Neta-Johnson and Neta of sixteenth order.

\begin{figure*}[!t]
\centering
\resizebox{\hsize}{!}{\includegraphics{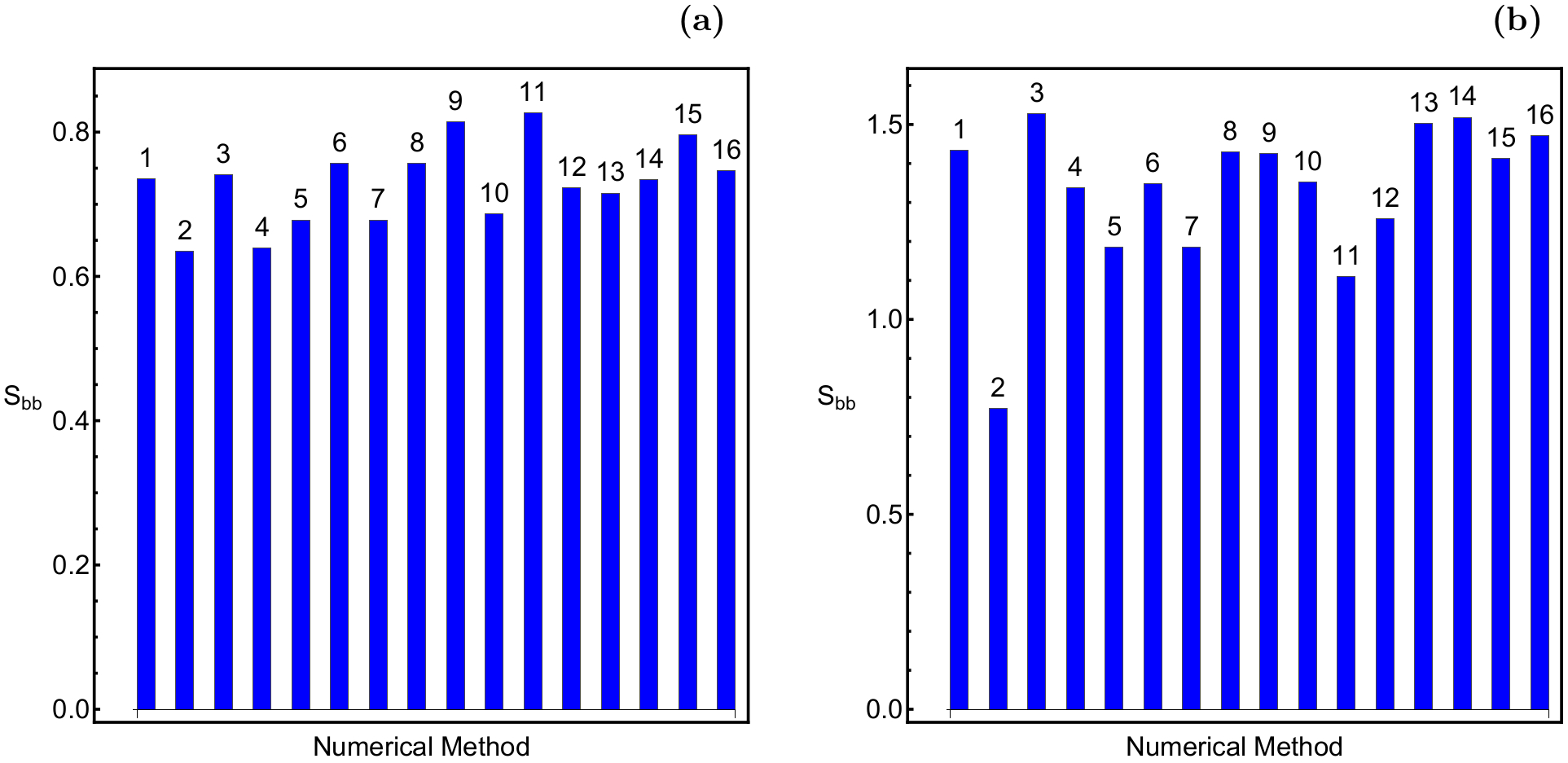}}
\caption{The distribution of the values of the boundary basin entropy $S_{bb}$ of the complex plane, for all the numerical methods for the case where (a-left): three roots exist and (b-right): nine roots exist. The numbers of each bar correspond to the numerical formulae, as they have been listed at the beginning of Section \ref{num}. (Color figure online).}
\label{bbe}
\end{figure*}

Despite the thorough study of the basins of convergence, their fractal boundaries are only visible in magnifications of local areas of the basin diagrams. In order to obtain quantitative results regarding the fractal nature of the basins we computed the boundary basin entropy \cite{DWGGS16}
\begin{equation}
S_{bb} = \frac{S}{N_b}
\label{sbb}
\end{equation}
where $N_b$ is the number of boxes containing more than one attractor.

According to the so-called ``log 2 criterion", the basin boundaries are not smooth (fractal) if $S_{bb} > \log 2$ (the converse is not true). In Fig. \ref{bbe} (a-b) we demonstrate the numerical values of the boundary basin entropy $S_{bb}$, for all the numerical methods. For both cases, regarding the total number of the roots, the lowest value of $S_{bb}$ corresponds to the Halley method. It is interesting to note that always (for all methods and for both cases) $S_{bb} > \log 2$, which implies that the basin boundaries on the complex plane are always fractal.

\section{Conclusions}
\label{conc}

In this work we used a large variety of numerical methods in order to reveal the basins of attraction on the complex plane, associated with the roots of simple polynomial equations. All the magnificent basin structures on the complex plane were identified by classifying dense grids of initial conditions, using the corresponding iterative schemes. In particular, we managed to determine how the geometry of the convergence structures changes as a function of the order of the applied numerical methods. Furthermore, all the correlations between the attracting domains and the corresponding distributions of the probability as well as the required number of iterations have been successfully established.

The following list contains the most important conclusions of our numerical analysis:
\begin{enumerate}
  \item In both cases, regarding the total number of the roots, and for all the studied numerical methods the corresponding basins of attraction extend to infinity.
  \item We detected several initial conditions for which the iterative formulae lead quickly to extremely large numbers. A small portion of these initial conditions do eventually converge to one of the roots, while for the vast majority of them the iterative procedure tend to infinity. This type of initial conditions is almost always present in the case where nine roots exist. On the other hand, where three roots are present this phenomenon is observed only for a few types of numerical methods.
  \item We identified a portion of initial conditions for which the numerical codes abort at the first step of the iterative procedure. These ill-behaved initial conditions are mainly present in the case where nine roots exits, while for the case where three roots are present they are encountered only in the Neta method of fourteenth order.
  \item In the case of three roots more than 98\% of the classified initial conditions per grid, for almost all the numerical methods, converge within the first 20 iterations. On the contrary, when nine roots exist the same amount of initial conditions require at least 30 iterations for converging, with the same desired accuracy.
  \item A clear relation, between the convergence speed of the iterative schemes and the order of the methods, has been observed, in the case of three roots. In particular, it was found that the convergence speed of the methods constantly decreases with increasing order. The same principle seems to hold also in the case with nine roots however, in that case we cannot draw safe conclusions.
  \item The most compact distributions of iterations were found to exist in the case of three roots, while in the case of nine roots the corresponding probability histograms exhibit extended tails. These extended tails, numerically expressed by the high values of the diversity $b$ as well as the differential entropy $h$, is the main reason of why the Laplace PDF fails to properly fit most of the probability histograms, in the case of nine roots.
  \item In both examined cases, regarding the total number of the roots, the numerical method displaying the lowest degree of fractality (or in other words the method with the most smooth basin boundaries) was the Halley method. On the other hand, the most chaotic method (with the most fractal basin boundaries) was the Maheshwari method, when three roots are present. In the case of nine roots, the numerical methods Chun-Neta, Neta-Johnson and Neta of sixteenth order were found to be almost equally chaotic.
\end{enumerate}

The general conclusion is that as the number of the roots (numerical attractors) increases the geometry of the complex plane becomes more complicated which leads to the following phenomena: (i) the amount of ill-behaved initial conditions increases, (ii) the required number of iterations increases, (iii) the degree of fractality, expressed through the basin entropy, increases.

A double precision numerical routine, written in standard \verb!FORTRAN 77! \cite{PTVF92}, was used for the classification of the initial conditions on the complex plane. Using a Quad-Core i7 2.4 GHz PC we needed about 4 minutes of CPU time, for performing the classification in each grid of initial conditions. Moreover, all the graphical illustration of the paper has been created using the latest version 11.3 of Mathematica$^{\circledR}$ \cite{W03}.

We hope that the presented numerical outcomes as well as the comparison of the several aspects of the numerical methods to be useful in the active field basins of attraction.

\section*{Acknowledgments}
\footnotesize

The authors would like to express their warmest thanks to the two anonymous referees for the careful reading of the manuscript and for all the apt suggestions and comments which allowed us to improve both the quality and the clarity of the paper.

\end{document}